\documentclass[12pt]{article}
\usepackage{times}
\usepackage{amsmath,amsxtra,amsthm,amssymb,makeidx,graphics,hyperref,graphicx}
\usepackage{latexsym,color}

\theoremstyle{plain}
\newtheorem{theorem}{Theorem}[section]
\newtheorem{lemma}[theorem]{Lemma}
\newtheorem{corollary}[theorem]{Corollary}

\theoremstyle{definition}
\newtheorem{remark}[theorem]{Remark}
\newtheorem{definition}[theorem]{Definition}
\newtheorem{example}[theorem]{Example}

\newtheorem{exercise}[theorem]{Exercise}

\theoremstyle{remark}

\newcommand{\ncom}{\newcommand}
  
\newcommand{\cA}{{\cal A}}

\newcommand{\cB}{{\cal B}}

\newcommand{\cE}{{\cal E}}
\newcommand{\cF}{{\cal F}}

\newcommand{\cI}{{\cal I}}

\newcommand{\cN}{{\cal N}}

\newcommand{\C}{{\mathbb C}}

\newcommand{\Fq}{{{\mathbb F}_q}}
\newcommand{\F}{{{\mathbb F}}}

\newcommand{\qb}[2]{{{ {{#1}\choose {#2}}_q }}}
\newcommand{\bin}[2]{{{ {#1}\choose {#2}}} }

\ncom{\n}[1]{{{ \|{#1}\|     }}}
\ncom{\ns}{\normalsize}
\ncom{\la}{\lambda}
\ncom{\bm}{\boldmath}
\ncom{\noi}{\noindent}
\ncom{\bq}{\begin{equation}}
\ncom{\eq}{\end{equation}}  
\ncom{\beqn}{\begin{eqnarray*}}
\ncom{\eeqn}{\end{eqnarray*}}  
\ncom{\ba}{\begin{array}}
\ncom{\ul}{\underline}   

\ncom{\ea}{\end{array}}
\ncom{\beq}{\begin{eqnarray}}
\ncom{\eeq}{\end{eqnarray}}  
\ncom{\nno}{\nonumber}
\ncom{\hs}{\mbox{\hspace{.25cm}}}
\ncom{\rar}{\rightarrow}
\ncom{\lrar}{\longrightarrow}
\ncom{\Rar}{\Rightarrow}
\ncom{\noin}{\noindent} 
\ncom{\bc}{\begin{center}}
\ncom{\ec}{\end{center}}  
\ncom{\sz}{\scriptsize}   
\ncom{\fpd}{\Phi(\pi^{'})}
\ncom{\fp}{\Phi(\pi) }
\ncom{\nk}{\left< \begin{array}{c}
                       n\\k \end{array} \right>}
\ncom{\nd}{1^{'},2^{'},\cdots,n^{'}}
\ncom{\de}{\bigtriangleup (F_{2n},\leq)}
\ncom{\del}{\bigtriangleup}
\ncom{\cov}{<\!\!\!\!\cdot }
\ncom{\bt}{\begin{theorem}}
\ncom{\bcon}{\begin{con}}
\ncom{\et}{\end{theorem}}
\ncom{\econ}{\end{con}}
\ncom{\bl}{\begin{lemma}}
\ncom{\el}{\end{lemma}}  
\ncom{\bco}{\begin{corollary}} 
\ncom{\ds}{\displaystyle}
\ncom{\eco}{\end{corollary}}   
\ncom{\bp}{\begin{pro}}  
\ncom{\ep}{\end{pro}}    
\ncom{\bex}{\begin{example}}
\ncom{\eex}{\end{example}}  
\ncom{\bexr}{\begin{exercise}}
\ncom{\eexr}{\end{exercise}}  
\ncom{\bprob}{\begin{problems}}
\ncom{\eprob}{\end{problems}}  
\ncom{\bd}{\begin{definition}}
\ncom{\ed}{\end{definition}}  
\ncom{\brm}{\begin{remark}}   
\ncom{\erm}{\end{remark}}     
\ncom{\bal}{\begin{Algorithm}}
\ncom{\eal}{\end{Algorithm}}  
\ncom{\ol}{\overline}
\ncom{\wh}{\widehat} 

\ncom{\pf}{\noi {\bf Proof.  }}
\ncom{\eprf}{\noi {$\Box$}}
\ncom{\be}{\begin{enumerate}} 
\ncom{\ee}{\end{enumerate}}   
\ncom{\seq}{\subseteq}

\ncom{\zr}{\bf\textcolor{red}}
\ncom{\zb}{\bf\textcolor{blue}}
\ncom{\zg}{\bf\textcolor{green}}
\ncom{\zm}{\bf\textcolor{magenta}}

\definecolor{gold}{rgb}{0.85,.66,0}
\definecolor{gb}{rgb}{0, .5,.5}
\definecolor{rb}{rgb}{0.5, 0,.5}
\definecolor{Pink}{rgb}{1,0.75, 0.8}
\ncom{\zgb}{\bf\textcolor{gb}}
\ncom{\zrb}{\bf\textcolor{rb}}

\makeatletter
\newcommand{\single}{\let\CS=\@currsize\renewcommand{\baselinestretch}{1.5}\tiny\CS}
\newcommand{\oneandahalfspacing}{\let\CS=\@currsize\renewcommand{\baselinestretch}{1.5}\tiny\CS}
\newcommand{\doublespacing}{\let\CS=\@currsize\renewcommand{\baselinestretch}{1.6}\tiny\CS}
\newcommand{\double}{\let\CS=\@currsize\renewcommand{\baselinestretch}{3}\tiny\CS}

\def\notarro{{{\hbox{{\hspace*{-.02in}}$\rightarrow$ } }
{\hbox{$\!\!\!\!\!\!\!$}}{\raise 0.2ex 
\hbox{$\scriptscriptstyle{/}$}}}\hspace{.06in}}
\def \*{^{\mbox{$*$}}}

\newcommand{\inp}[2]{\langle {#1} ,\,{#2} \rangle}

\newcommand\bnota{\begin{nota} }                         
\newcommand\enota{\end{nota} }                           

\newcommand{\beano}{\begin{eqnarray*}}
\newcommand{\eeano}{\end{eqnarray*}}

\begin{document}

\title{\bf{{{A $q$-analog of the adjacency matrix of the $n$-cube}}}}

\author{{\bf 
{Subhajit Ghosh}} \\
{\em  {Department of Mathematics}}\\
{\em  {Bar-Ilan University}}\\
{\em  {Ramat-Gan, 5290002 ISRAEL}}\\
{\bf  \texttt{gsubhajit@alum.iisc.ac.in}}\\ $\;$\\
{\bf {Murali K. Srinivasan}} \\
{\em  {Department of Mathematics}}\\
{\em  {Indian Institute of Technology, Bombay}}\\
{\em  {Powai, Mumbai 400076, INDIA}}\\
{\bf  \texttt{murali.k.srinivasan@gmail.com}}}
 
\date{}
\maketitle
\begin{center}
{\small{\em To the memory of Chandan and
beloved Lawson}}
\end{center}

\begin{abstract} 

Let $q$ be a prime power and define $(n)_q = 1+q+q^2+\cdots +q^{n-1}$, 
for a nonnegative integer $n$. Let $B_q(n)$ denote the set 
of all subspaces of  $\F_q^n$, 
the $n$-dimensional $\Fq$-vector space of all column vectors 
with $n$ components. 

Define a $B_q(n)\times B_q(n)$ complex matrix $M_q(n)$  
with entries given by
	\beqn
	M_q(n)(X,Y) &=& \left\{ \ba{ll}  
	1 & \mbox{if }Y\seq X, \dim(Y)= \dim(X)-1,\\
	q^k & \mbox{if }X\seq Y, \dim(Y)=k+1, \dim(X)=k,\\
	0 & \mbox{otherwise.} \ea\right.
	\eeqn
We think of $M_q(n)$ as a $q$-analog of the 
adjacency matrix of the $n$-cube. 
We show that the eigenvalues of $M_q(n)$ are
$$(n-k)_q - (k)_q 
\mbox{ with multiplicity }  \qb{n}{k},\;\;k=0,1,\ldots ,n,$$
and we write down an explicit canonical eigenbasis of $M_q(n)$.
We give a weighted count of
the number of rooted spanning trees in the $q$-analog of the $n$-cube.

\end{abstract}
     
\noi

\section{{Introduction}}

One aspect of algebraic combinatorics is the study of
eigenvalues and eigenvectors of certain matrices 
associated with posets and graphs. Among the most basic such examples is 
the adjacency matrix of the $n$-cube, which 
has an elegant spectral theory and arises in a variety of applications
(see {\bf\cite{cst,st}}). This paper defines a 
$q$-analog of this matrix, studies its spectral theory, and gives
an application to weighted counting of rooted spanning trees in the $q$-analog
of the $n$-cube.

Let $q$ be a prime power and define $(n)_q = 1+q+q^2+\cdots +q^{n-1}$, 
for a nonnegative integer $n$. Let $B_q(n)$ denote the set 
of all subspaces of  $\F_q^n$, 
the $n$-dimensional $\Fq$-vector space of all column vectors 
with $n$ components. 
The set of $k$-dimensional subspaces in $B_q(n)$ is denoted $B_q(n,k)$
and its cardinality is the {\em $q$-binomial coefficient $\qb{n}{k}$}.  
The {\em Galois number}
$$G_q(n)=\sum_{k=0}^n \qb{n}{k}$$
is the total number of subspaces in $B_q(n)$. The set $B_q(n)$ has 
the structure of a graded poset of rank $n$, under inclusion.

Recall the definition of the adjacency matrix $M(n)$ of the $n$-cube: let
$B(n)$ denote the set of all subsets of $\{1,2,\ldots ,n\}$. The 
rows and columns of $M(n)$ are indexed by elements of $B(n)$,
with the entry in row $S$, column $T$ equal to 1 if $|(S\setminus T)\cup 
(T\setminus S)|=1$
and equal to 0 otherwise.
Define a $B_q(n)\times B_q(n)$ complex matrix $M_q(n)$  with entries given by
	\begin{eqnarray}\label{eq:mt1}
	M_q(n)(X,Y) &=& \left\{ \ba{ll}  
	1 & \mbox{if }Y\seq X, \dim(Y)= \dim(X)-1,\\
	q^k & \mbox{if }X\seq Y, \dim(Y)=k+1, \dim(X)=k,\\
	0 & \mbox{otherwise.} \ea\right.
	\end{eqnarray}
We think of $M_q(n)$ as {\em a $q$-analog of the adjacency matrix of the $n$-cube}. 
Note that $M(n)$  
is symmetric, has entries in $\{0,1\}$, and has all
row sums equal to $n$. The significance of the definition above for the $q$-analog
comes from the fact that though
$M_q(n)$ lacks the first two properties
it does have all row sums equal.
Indeed, let $X\in B_q(n,k)$. Then the number of subspaces covering
$X$ is $(n-k)_q$ and the number of subspaces covered by $X$ is $(k)_q$ and so
the sum of the entries of row $X$ of
$M_q(n)$ is $q^k(n-k)_q + (k)_q = (n)_q$. 

A scaling (i.e., a diagonal similarity) of $M_q(n)$ is symmetric. Let $D_q(n)$ be the $B_q(n)\times B_q(n)$
diagonal matrix with diagonal entry in row $X$, column $X$ 
given by $\sqrt{q^{\bin{k}{2}}}$, where $k=\dim(X)$. Then for $X\in B_q(n,k)$,
$Y\in B_q(n,r)$ the entry in row $X$, column $Y$ of $D_q(n)M_q(n)D_q(n)^{-1}$
is given by
        \beqn
        \lefteqn{\sqrt{q^{\bin{k}{2}}}\, M_q(X,Y)\, \sqrt{q^{-\bin{r}{2}}}}\\
        &=&\left\{\ba{ll}
        \sqrt{q^{\bin{k}{2}}}\, q^k\, \sqrt{q^{-\bin{k+1}{2}}} =
        \sqrt{q^k} & \mbox{if $X\seq Y$ and $r=k+1$,}\\   
        \sqrt{q^{\bin{k}{2}}} \sqrt{q^{-\bin{k-1}{2}}} =  
        \sqrt{q^{r}} & \mbox{if $Y\seq X$ and $r=k-1$,}\\ 
        0 & \mbox{otherwise,} 
        \ea
        \right.\\
        &=&\left\{\ba{ll}
        \sqrt{q^{\min\{\dim(X),\dim(Y)\}}}&
\mbox{if $X\seq Y$ or $Y\seq X$, and $|\dim(X)-\dim(Y)|=1$,} \\
        0 & \mbox{otherwise,} 
        \ea
        \right.
        \eeqn  
yielding a symmetric matrix. It follows that $M_q(n)$
is diagonalizable and that its eigenvalues 
are real. In fact they are integral and the eigenvalue-multiplicity
pairs of $M_q(n)$ are a $q$-analog of those for $M(n)$; the eigenvalues of
$M(n)$ are $n-2k=(n-k)-(k)$ with multiplicity 
$\binom{n}{k}$, $k=0,1,\ldots ,n$ (see{\bf\cite{cst,st}}).

\bt \label{mt1}
The eigenvalues of the matrix $M_q(n)$ are
$$(n-k)_q - (k)_q 
\mbox{ with multiplicity }  \qb{n}{k},\;\;k=0,1,\ldots ,n.$$
\et

We give two proofs of Theorem \ref{mt1} in this paper. 
In Section 2 we use a result of Terwilliger {\bf\cite{t1}} on the up operator
on subspaces to show that Theorem \ref{mt1} reduces
to showing  that the eigenvalues of $K_q(n)$, 
a certain $(n+1)\times (n+1)$ tridiagonal matrix, are 
$(n-k)_q - (k)_q,\;0\leq k \leq n$. The matrix
$K_q(n)$ is a $q$-analog 
of the famous tridiagonal matrix $K(n)$ of 
Mark Kac {\bf\cite{a,k,tt}} (both these
matrices are defined in Section 2). Now $K_q(n)$ occurs in 
Terwilliger's classification
of Leonard pairs ({\bf\cite{t2, t4}}) and as such its eigenvalues/eigenvectors
were known (see, for example, Lemma 4.20 in {\bf\cite{t3}}), completing
the first proof of Theorem \ref{mt1}.

Another natural problem is to write down eigenvectors of $M_q(n)$.
For $X\in B_q(n)$ with $\dim(X)=k$ define
$$\pi(X) = \frac{q^{\bin{k}{2}}}{P_q(n)},
$$
where $P_q(n) =\prod_{k=0}^{n-1} (1+q^k)$.
We have
$$\sum_{X\in B_q(n)}\pi(X)\;=\;\frac{\sum_{k=0}^n q^{\bin{k}{2}}\qb{n}{k}}
{P_q(n)}\;=\;1$$
where the second equality follows by
the $q$-binomial theorem (so  $\pi$ is a probability vector on $B_q(n)$).

Define an inner product on the (complex) vector space of column vectors
with components indexed by $B_q(n)$ as follows: given vectors $u,v$ define
\beq \label{inp}
\inp{u}{v}_\pi =\sum_{X\in B_q(n)}\ol{u(X)}v(X)\pi(X).
\eeq
Since $P_q(n)$ is independent of $k$,
the argument showing that $D_q(n)M_q(n)D_q(n)^{-1}$
is symmetric shows that $M_q(n)$ is
self-adjoint with respect
to the inner product (\ref{inp}).

Recall that a classical result exhibits an explicit orthogonal eigenbasis
of $M(n)$ (under the standard inner product), see {\bf\cite{cst,st}}. Up to
scalars, this basis is canonical in the sense that no choices are involved
in writing it down. 
In Section 4 we extend this method to $M_q(n)$. The main idea is to use
a linear algebraic interpretation of the Goldman-Rota recurrence for the
Galois numbers (see {\bf\cite{gr,kc,ku}}) that was given in {\bf\cite{s3}} 
(and summarized in Section 3).

\bt \label{evectors}
There is an inductive procedure to write down a canonical eigenbasis 
of $M_q(n)$, orthogonal with respect to the inner product (\ref{inp}). 
\et
In the course of writing down the eigenvectors we also evaluate the matrix
eigenvector product, thereby 
giving an alternate proof of Theorem \ref{mt1}.

We originally arrived at the matrix $M_q(n)$  
through a reversible Markov chain with state space
$B_q(n)$, transition matrix $\frac{1}{(n)_q}M_q(n)$, and stationary distribution
$\pi$ (see {\bf\cite{gs}}). 
Since the spectral
theory of $M_q(n)$, as a $q$-analog of the spectral theory of $M(n)$, 
is of independent interest we are presenting it separately
in this paper. 
Another application concerns weighted enumeration of 
rooted spanning trees in the $q$-analog of the $n$-cube. 

Let us first recall 
the remarkable product formula (see Example 9.12 in {\bf\cite{st}}) for the
number of rooted spanning trees in the n-cube. 
Given a multigraph $G$ 
without loops and a vertex $v$ in $G$, let $\tau(G,v)$ denote the number
of rooted spanning trees in $G$ with root $v$ and let $\tau(G)$ denote the total
number of rooted spanning trees in $G$. Let $C(n)$ denote the $n$-cube (i.e.,
the Hasse diagram of the poset $B(n)$ treated as a graph). The Laplacian
eigenvalues of $C(n)$ are well known to be $2k,\;k=0,1,\ldots ,n$ with multiplicity
$\bin{n}{k}$. It thus follows from the Matrix-Tree theorem that
\beq \label{hc}
\tau(C(n))&=&\prod_{k=1}^n (2k)^{\bin{n}{k}}.
\eeq
The {\em $q$-analog $C_q(n)$} of $C(n)$ is defined to be the Hasse diagram
of $B_q(n)$ treated as a graph (note that $M_q(n)$ is not the adjacency
matrix of the grapn $C_q(n)$). 
The eigenvalues of the Laplacian of $C_q(n)$ are not known
(note that $C_q(n)$ is not regular) but it was shown in {\bf\cite{s2}} that
the product of the nonzero eigenvalues of the Laplacian is more tractable and
led to a product formula for the number of spanning trees 
of $C_q(n)$, although the individual terms in this product are 
not explicilty given but only as a positive combinatorial sum.
Here we obtain an explicit $q$-analog of (\ref{hc}) by a weighted
count of the rooted spanning trees of $C_q(n)$. 

Let $F\in \cF_q(n)$, the set of all 
rooted spanning trees of $C_q(n)$. Orient every edge
of $F$ by pointing it towards the root. Let $e=(X,Y)$ be an oriented edge
of $F$.  We say $e$ is {\em{spin up}} if $\dim(Y)=\dim(X)+1$ and is {\em{spin down}}
if $\dim(Y)=\dim(X)-1$. The {\em{weight}} of $F$ is defined by
\beqn
w(F)&=&\sum_{(X,Y)}\dim(X),
\eeqn
where the sum is over all spin up oriented edges of $F$. In Section 5 we prove
the following result 
(the proof can be read at this point, assuming Theorem \ref{mt1}).
\bt \label{qct} We have
$$\sum_{F\in\cF_q(n)}q^{w(F)} = \prod_{k=1}^n ((1+q^{n-k})(k)_q)^{\qb{n}{k}}.
$$
\et

\section{{Eigenvalues of $M_q(n)$}}

The spectral theory of $M_q(n)$ goes hand in hand with that of
a $q$-analog of the Kac matrix.
Recall that the Kac matrix is a $(n+1)\times (n+1)$ tridiagonal matrix $K(n)$
with diagonal $(0,0,\ldots ,0)$, subdiagonal $(1,2,\ldots ,n)$ 
and superdiagonal $(n,n-1,\ldots ,1)$: 
\[K(n) = \left[\ba{cccccc}
0&n&&&&\\
1&0&n-1&&&\\
&2&0&n-2&&\\
&&\ddots&\ddots&\ddots&\\
&&&n-1&0&1\\
&&&&n&0
\ea\right].\]
The eigenvalues of $K(n)$ are $n-2k, k=0,1,\ldots ,n$, and
its eigenvectors have been written down 
(see {\bf\cite{a,cst,ek,k,tt}}).

We define the $q$-analog of the Kac matrix to be 
the $(n+1)\times (n+1)$ tridiagonal matrix $K_q(n)$
with diagonal $(0,0,\ldots ,0)$, subdiagonal $((1)_q,(2)_q,\ldots ,(n)_q)$ and superdiagonal $((n)_q,q(n-1)_q,\ldots ,q^{n-1}(1)_q)$: 
\[K_q(n) = \left[\ba{cccccc}
0&(n)_q&&&&\\
(1)_q&0&q(n-1)_q&&&\\
&(2)_q&0&q^2(n-2)_q&&\\
&&\ddots&\ddots&\ddots&\\
&&&(n-1)_q&0&q^{n-1}(1)_q\\
&&&&(n)_q&0
\ea\right].\]
More formally, let us index the rows and columns of $K_q(n)$ by the set
$\{0,1,2,\ldots ,n\}$.
If $c_0,c_1,\ldots ,c_n$ denote column vectors
in $\F_q^{n+1}$ with $c_i$ having a 1 in the component 
indexed by $i$ and 0's
elsewhere (and we set $c_{-1}=c_{n+1}=0$) 
then, for  $0\leq k \leq n$, column $k$ of $K_q(n)$
is 
\beq \label{ctm}
&(k+1)_q\;c_{k+1} + q^{k-1} (n-k+1)_q\;c_{k-1}.&
\eeq

\bt \label{mt2}
(i) The eigenvalues of $K_q(n)$ are 
$$(n-k)_q - (k)_q,\;k=0,1, \ldots ,n.$$

(ii) For $0\leq k \leq n$, there is a right eigenvector of $K_q(n)$
corresponding to the eigenvalue $(n-k)_q - (k)_q$ whose component 
$i,\;0\leq i \leq n$ is given by 
$$ _3\phi_2 \left( \ba{c} q^{-i}, q^{-k}, -q^{k-n}\\0, q^{-n} \ea
 \; \vline \; q, q            
\right)$$
where the basic hypergeometric series notation is from {\bf\cite{gr1}}. 
\et
\pf 
This is Lemma 4.20 in Terwilliger {\bf\cite{t3}} (replace $b$ by $q$,
$d$ by $n$, and $j$ by $n-k$).
Part (i) also follows from 
Theorem 2 in Johnson {\bf\cite{j}} (by 
taking $a=z=1$ and $h=k=0$). \eprf

We now relate the spectral theories of $M_q(n)$ and $K_q(n)$.
The proper framework for studying this
are the up (and down) operators on the poset of subspaces.

For a finite set $S$, we denote the complex vector space with $S$
as basis by $\C[S]$. 
We denote by $r$ the rank function (given by dimension) of the graded
poset $B_q(n)$.
Then we have (vector space direct sum)
$$\C[B_q(n)]=\C[B_q(n,0)]\oplus 
\C[B_q(n,1)]\oplus \cdots \oplus\C[B_q(n,n)].$$

An element $v\in \C[B_q(n)]$ is {\em homogeneous} if 
$v\in \C[B_q(n,i)]$ for some $i$,
and if $v\not= 0$, we extend the notion of rank to nonzero 
homogeneous elements by writing $r(v)=i$. For $0\leq k \leq n$,
the {\em $k^{th}$ up operator}  $U_{n,k}:\C[B_q(n)]\rar \C[B_q(n)]$ 
is defined, for $X\in B_q(n)$, by $U_{n,k}(X)=0$ if $\dim(X)\not=k$ and
$U_{n,k}(X)= \sum_{Y} Y$,
where the sum is over all $Y\in B_q(n)$ covering $X$, if $\dim(X)=k$.
Similarly we define the {\em $k^{th}$ down operator} 
$D_{n,k}:\C[B_q(n)]\rar \C[B_q(n)]$ (we have $U_{n,n}=D_{n,0}=0$).
Set $U_n=U_{n,0}+U_{n,1}+\cdots +U_{n,n}$ 
and $D_n=D_{n,0}+D_{n,1}+\cdots +D_{n,n}$, called, respectively,
the {\em up and down} operators on $\C[B_q(n)]$.
For a finite vector space $X$ over $\Fq$ we denote by
$B_q(X)$ the set of all subspaces of $X$ 
and we denote by $U_X$ (respectively, $D_X$) the up
operator (respectively, down operator) on $\C[B_q(X)]$.

If we think of the elements of $\C[B_q(n)]$ as column vectors 
with components
indexed by the {\em standard basis} elements $B_q(n)$ then 
$M_q(n)$ is the matrix of the operator
$$\cA_q(n) = U_n +\sum_{k=0}^n q^{k-1}D_{n,k}$$
with respect to the basis $B_q(n)$. 

For $0\leq k \leq n$, define $s_k\in \C[B_q(n,k)]$ by
$$ s_k = \sum_{X\in B_q(n,k)} X,$$
and define $R_q(n)$ to be the subspace of $\C[B_q(n)]$ spanned by $s_0,s_1,\ldots
, s_n$. Elements of $R_q(n)$ are called {\em radial vectors}. Clearly, $R_q(n)$
is closed under $\cA_q(n)$ and $\dim(R_q(n))=n+1$. We have
\beqn
\cA_q(n)(s_k)= (k+1)_q\;s_{k+1} + q^{k-1}(n-k+1)_q\;s_{k-1},\;\;0\leq k \leq n.
\eeqn
It follows from (\ref{ctm})
that the matrix of $\cA_q(n): R_q(n)\rar R_q(n)$ with respect to the basis
$\{s_0,\ldots ,s_n\}$ is $K_q(n)$.

A {\em symmetric chain} in $\C[B_q(n)]$ is a sequence
\beq \label{gjc}
&s=(v_k,\ldots ,v_{n-k}),\;\;\;k\leq n/2,&
\eeq
of nonzero homogeneous elements of $\C[B_q(n)]$
such that 
\begin{itemize}
\item $r(v_i)=i$ for $i=k,\ldots ,n-k$.

\item $U_n(v_i)$ is a nonzero scalar multiple of $v_{i+1}$, for
$i=k,\ldots , n-k-1$ and $U_n(v_{n-k})=0$.

\item $D_n(v_{i+1})$ is a nonzero scalar multiple of $v_i$ for
$i=k,\ldots , n-k-1$ and $D_n(v_k)=0$.
\end{itemize}
Note that the
elements of the sequence $s$ are linearly independent, being nonzero and of
different ranks. 
We say that $s$ {\em
starts} at rank $k$ and {\em ends} at rank $n-k$. Note that the subspace spanned
by the elements of $s$ is closed under $U_n,D_n$ and also $\cA_q(n)$.

The following result was proved (in an equivalent form) 
by Terwilliger {\bf\cite{t1}} (see Item 5 of Theorem 3.3 on top
of page 208). For a proof
using Proctor's $\mathfrak{sl}(2,\C)$ technique {\bf\cite{p}} see Theorem 2.1 
in {\bf\cite{s2}} (where also the result is stated differently
but in an equivalent form to that given below).

\bt \label{sjb}
There exists a basis $T_q(n)$ of $\C[B_q(n)]$ such that
\begin{enumerate}
\item $T_q(n)$ is a disjoint union of symmetric chains in $\C[B_q(n)]$.

\item Let $0\leq k \leq n$  and let $(v_k,\ldots ,v_{n-k})$ be any
symmetric chain in $T_q(n)$ starting at rank $k$ and ending at rank $n-k$.
Then
\beqn
U_n(v_u)&=& q^k(u+1-k)_q\;v_{u+1},\;k\leq u < n-k.\\
D_n(v_{u+1})&=& (n-k-u)_q\;v_u,\;k\leq u < n-k.
\eeqn
\end{enumerate}
\et

We now give the 

\noi {\bf{Proof of Theorem \ref{mt1}}} Observe the following.

(i) The number of symmetric chains in $T_q(n)$ starting at rank $k$ and ending at
rank $n-k$, for $0\leq k \leq n/2$, is $\qb{n}{k} - \qb{n}{k-1}$.

(ii) Let $s=(v_k,\ldots ,v_{n-k})$ be a symmetric chain in $T_q(n)$ starting
at rank $k$, where
$0\leq k \leq n/2$. Then the subspace spanned by $\{v_k,\ldots ,v_{n-k}\}$ 
is closed under $\cA_q(n)$ and the matrix of $\cA_q(n)$ 
with respect to the basis $s$ is $q^k K_q(n-2k)$.

(iii) By Theorem \ref{mt2} the eigenvalues of $q^k K_q(n-2k)$ are
\beqn
\lefteqn{ q^k((n-2k-i)_q - (i)_q),\;i=0,1,\ldots ,n-2k}\\ 
 &=& ((n-i)_q - (i)_q),\;i=k,\ldots ,n-k.
\eeqn

(iv) It now follows from items (i), (ii), (iii) above that the eigenvalues
of $\cA_q(n)$ are
\beqn
(n-j)_q - (j)_q,\;j=0,\ldots ,n,
\eeqn
with respective multiplicities
\beqn \sum_{i=0}^{\min\{j,n-j\}} \qb{n}{i} - \qb{n}{i-1} = \qb{n}{j}.
\eeqn
That completes the proof. \eprf

\section{{A decomposition of $\C[B_q(n)]$}}

In this and the next section we give proofs of Theorems \ref{mt1}
and \ref{evectors} by inductively writing down 
an eigenbasis for the operator $\cA_q(n)$.
This is based on
a direct sum decomposition of the vector space $\C[B_q(n)]$ that was worked out
in the paper {\bf\cite{s2}}. This decomposition yields a linear 
algebraic interpretation 
of the Goldman-Rota recurrence  
for the Galois numbers and is
of independent interest.
In {\bf\cite{s2}} it was used to inductively
write down an explicit eigenbasis for the Bose-Mesner algebra 
of the Grassmann scheme. Here
we recall the relevant definitions and results from {\bf\cite{s2}}. All the 
omitted proofs may be found in Section 2 of {\bf\cite{s2}}.

The Goldman-Rota identity {\bf\cite{gr,kc,ku}} is the recursion
\beq \label{gri}
G_q(n+1) &=& 2G_q(n) + (q^n - 1)G_q(n-1),\;n\geq 1,\eeq
with $G_q(0)=1,\;G_q(1)=2,$
or, in terms of the $q$-binomial coefficient,
\beq \label{nsw}
\qb{n+1}{k} &=& \qb{n}{k} + \qb{n}{k-1} + (q^n - 1)\qb{n-1}{k-1},\;n,k\geq
1,\eeq
with $\qb{0}{k}=\delta(0,k)$ and $\qb{n}{0}=1$.
Note that (\ref{gri}) follows by summing (\ref{nsw}) over $k$.

We shall now give a linear algebraic interpretation to (\ref{gri}). Denote
the standard basis vectors of $\F^n_q$ by the column vectors $e_1,\ldots ,e_n$.  
We identify $\F_q^k$, for $k<n$, with the subspace of $\F_q^n$ consisting  of
all vectors with the last $n-k$ components zero. 
So $B_q(k)$ consists of all
subspaces of $\F_q^n$ contained in the subspace spanned by $e_1,\ldots ,e_k$.

Define $A_q(n)$ to be the
collection of all subspaces in $B_q(n)$ not contained in the hyperplane
$\F_q^{n-1}$, i.e.,
$$A_q(n) = B_q(n) - B_q(n-1) = \{ X\in B_q(n) : X\not\subseteq \F_q^{n-1} 
\},\;n\geq 1.$$
For $1\leq k \leq n$, let $A_q(n,k)$ denote the set of all subspaces in
$A_q(n)$ with dimension $k$. We consider $A_q(n)$ 
as an induced subposet of $B_q(n)$.

We have a direct sum decomposition
\beq \label{bod}
\C[B_q(n+1)] = \C[B_q(n)] \oplus \C[A_q(n+1)].
\eeq
We shall now give a further decomposition of $\C[A_q(n+1)]$.

Let $H(n+1,\F_q)$ denote the subgroup of $GL(n+1,\F_q)$ consisting of all
matrices of the form
$$ \left[ \ba{cc}
           I & \ba{c} a_1 \\  \cdot \\ \cdot \\  a_n \ea \\
           0 \cdots 0 & 1
          \ea 
   \right],   
$$
where $I$ is the $n\times n$ identity matrix.

The additive abelian group $\F_q^n$ is isomorphic to $H(n+1,\F_q)$ via
$\phi : \F_q^n \rar H(n+1,\F_q)$ given by
$$\phi \left(\left[ \ba{c} a_1 \\ \cdot \\ \cdot \\a_n \ea \right]\right)
\rar
\left[ \ba{cc}
           I & \ba{c} a_1 \\  \cdot \\ \cdot \\  a_n \ea \\
           0 \cdots 0 & 1
          \ea 
   \right].   
$$   

There is a natural (left) action
of $H(n+1,\F_q)$ on $A_q(n+1)$ and $A_q(n+1,k)$. 

Let $\cI_q(n)$  
denote the set of all distinct irreducible characters (all of
degree $1$) of $H(n+1,\F_q)$
and let $\cN_q(n)$  
denote the set of all distinct nontrivial irreducible characters
of $H(n+1,\F_q)$.

For $\chi\in \cI_q(n)$,
let $W(\chi)$ (respectively, $W(\chi,k)$) denote the isotypical component of
$\C[A_q(n+1)]$ (respectively, $\C[A_q(n+1,k)]$) corresponding to the 
irreducible representation of $H(n+1,\F_q)$ with character $\chi$.    
When $\chi$ is the trivial character we denote $W(\chi)$ (respectively,
$W(\chi,k)$) by  $W(0)$ (respectively, $W(0,k)$).
We have the following decompositions, (note that $W(\chi,n+1)$, 
for $\chi\in \cN_q(n)$, is the zero module). 
\beq
W(0) &=& W(0,1) \oplus \cdots \oplus W(0,n+1), \\
W(\chi) &=& W(\chi,1) \oplus \cdots \oplus W(\chi,n),\;\;\;\chi \in \cN_q(n), \\
\C[A_q(n+1)] &=& W(0) \oplus \left(\oplus_{\chi \in \cN_q(n)} W(\chi)\right).
\eeq
Now $GL(n+1,\F_q)$ acts on $B_q(n+1)$ and
$U_{n+1}$ is $GL(n+1,\F_q)$-linear (and hence $H(n+1,\F_q)$-linear). 
Also, $\C[A_q(n+1)]$ is clearly closed under $U_{n+1}$. Thus 
\beq \label{d0}
W(0),\;\;W(\chi),\chi\in \cN_q(n) \mbox{ are  $U_{n+1}$-closed}.
\eeq

Define an equivalence 
relation $\sim$ on $A_q(n)$ by $X\sim Y$ iff 
$X\cap \,\F_q^{n-1}= Y\cap \,\F_q^{n-1}$.
Denote the
equivalence class of $X\in A_q(n)$ by $[X]$. 
For a subspace $X\in B_q(n-1)$, define $\wh{X}$ to be the subspace in $A_q(n)$
spanned by $X$ and $e_n$.

\bl \label{el}
Let $X,Y \in A_q(n)$ and $Z,T\in B_q(n-1)$. Then

\noi (i) $\mbox{dim}\,(X\cap \,\F_q^{n-1})=\mbox{dim}\,X - 1$ and 
$\wh{X\cap \,\F_q^{n-1}}\in [X]$. 

\noi (ii) $Z\leq T$ iff $\wh{Z} \leq \wh{T}$.

\noi (iii) $Y$ covers $X$ iff 

(a) $Y\cap \,\F_q^{n-1}$ covers $X\cap \,\F_q^{n-1}$ and

(b) $Y = \mbox{ span}\,((Y\cap \F_q^{n-1})\cup \,\{v\})$ for any
$v\in X - \F_q^{n-1}$. 

\noi (iv) The number of subspaces $Z'\in A_q(n)$ with
$Z'\cap \,\F_q^{n-1} = Z$ is $q^l$, where $l= n-\mbox{dim}\,Z - 1$.
Thus, $|[X]|=q^{n-k}$, where $k=\mbox{ dim}\,X$.\;\;\mbox{\eprf}
\el

For $X\in A_q(n+1)$, let
$G_X \subseteq H(n+1,\F_q)$ denote the stabilizer of $X$.

\bl \label{orl}
Let $X,Y\in A_q(n+1)$. Then

\noi (i) The orbit of $X$ under the action of $H(n+1,\F_q)$ is $[X]$.

\noi (ii) Suppose $Y$ covers $X$. Then the bipartite graph of the covering
relations between $[Y]$ and $[X]$ is regular with degrees $q$ (on the
$[Y]$ side) and $1$ (on the $[X]$ side). 

\noi (iii) Suppose $X\subseteq Y$. Then $G_X \subseteq G_Y$.\;\;\eprf
\el

Consider $\C[B_q(n+1)]$. For $X\in B_q(n)$ define
\beqn\theta_n(X)= \sum_Y Y,
\eeqn
where the sum is over all $Y\in A_q(n+1)$ covering $X$. Equivalently, the sum
is over all $Y\in A_q(n+1)$ with $Y\cap\,\F_q^n = X$, i.e., $Y\in [\wh{X}]$. It
follows from Lemma \ref{orl}(i) that 
\beq \label{d1}
\theta_n : \C[B_q(n)] \rar W(0)
\eeq
is a linear isomorphism.

We have the decomposition
\beq \label{d2}
\C[B_q(n+1)] &=& (\C[B_q(n)]\oplus W(0)) \oplus 
\left(\oplus_{\chi \in \cN_q(n)} W(\chi)\right),
\eeq
where, by (\ref{d0}) and (\ref{d1}), 
\beq \label{d3}
\C[B_q(n)]\oplus W(0)\mbox{ is  $U_{n+1}$-closed}.
\eeq

Let $\psi_k$ (respectively, $\psi$) denote the character of the
permutation
representation of $H(n+1,\F_q)$ on $\C[A_q(n+1,k)]$ (respectively,
$\C[A_q(n+1)]$) corresponding to the left action. Clearly $\psi = \sum_{k=1}^{n+1} \psi_k$.
Below $[,]$
denotes character inner product and the $q$-binomial coefficient $\qb{n}{k}$ is taken to
be zero when $n$ or $k$ is $< 0$.

\bt \label{grl}
\noi (i) Let $\chi \in \cI_q(n)$ be the trivial character. Then
$[\chi , \psi_k ]= \qb{n}{k-1},\;1\leq k \leq n+1.$

\noi (ii) Let $\chi \in \cN_q(n)$. Then
$[\chi , \psi_k ]= \qb{n-1}{k-1},\;1\leq k \leq n+1.$  \eprf
\et
Using Theorem \ref{grl}(ii) we see that
\beq \label{grd}
\dim(W(\chi))&=&\sum_{k=1}^{n+1} \qb{n-1}{k-1}=\;G_q(n-1),\;\;\chi\in\cN_q(n)
\eeq 
Now, by taking dimensions on both sides of (\ref{d2}) 
and using (\ref{d1}), (\ref{grd}) we get the Goldman-Rota identity (\ref{gri}).
More generally, by restricting to dimension $k$ on both sides of
(\ref{d2}), we get the identity (\ref{nsw}).

For $\chi \in \cI_q(n)$, define the following element of the group algebra of
$H(n+1,\F_q)$:
$$p(\chi)=\sum_g \ol{\chi(g)}\,g,$$
where the sum is over all $g\in H(n+1,\F_q)$. For $1\leq k \leq n+1$, the map
\beq \label{proj}
&p(\chi) : \C[A_q(n+1,k)] \rar \C[A_q(n+1,k)],&
\eeq
given by $v\mapsto \sum_{g\in H(n+1,\F_q)} \ol{\chi(g)}\,gv$,
is a nonzero multiple of the $H(n+1,\F_q)$-linear projection onto
$W(\chi,k)$. Similarly for $p(\chi) : \C[A_q(n+1)] \rar \C[A_q(n+1)]$.

For future reference we record the following
observation:
\beq \label{ds}
p(\chi)(\wh{Y}) \mbox{ and } p(\chi)(\wh{Z}) \mbox{ have disjoint supports, for }
Y\not= Z \in B_q(n).
\eeq

\bl \label{cl} Let $X\in A_q(n+1)$ and $\chi\in \cI_q(n)$. Then $p(\chi)(X)\not=
0$ iff $\chi : G_X \rar \C^*$ is the trivial character of $G_X$. \eprf
\el

\bt \label{f}

(i) Let $\chi\in \cI_q(n),\;X,Y\in A_q(n+1)$ with $X=hY$ for
some $h\in H(n+1,\F_q)$. Then
$$p(\chi)(X) = \ol{\chi(h^{-1})}\,p(\chi)(Y).$$

\noi (ii) Let $\chi\in \cI_q(n)$. Then
$\{ p(\chi)(\wh{X}) : X\in B_q(n,k-1) \mbox{ with } p(\chi)(\wh{X})
\not= 0\}$ is a basis of $W(\chi,k)$, $1\leq k \leq n+1$.

\noi (iii) Let $\chi\in \cI_q(n)$ and let $X,Y\in B_q(n)$ with $X$ covering $Y$.
$$p(\chi)(\wh{X})\not= 0 \mbox{ implies } p(\chi)(\wh{Y})\not= 0.\;\;\mbox{\eprf}$$

\et

Let $\chi \in \cN_q(n)$. 
By Theorem \ref{grl}(ii) we have $\dim(W(\chi,n))=1$.
It thus follows by Theorem \ref{f}(ii) and (\ref{ds}) above 
that there is a unique element $X(\chi)\in B_q(n,n-1)$
such that $p(\chi)(\wh{X(\chi)})\not=0$. Moreover,

\bl \label{cu}
Let $Y\in B_q(n,n-1)$. Then
$$|\{\chi\in\cN_q(n)\;|\; X(\chi)=Y\}| = q-1.\;\;\mbox{ \eprf}$$
\el 

\section{{Eigenvectors of $\cA_q(n)$}}\label{sec:Eig_vec}

In this section we inductively write down an
eigenbasis of $\cA_q(n)$ (in the process giving an alternate proof of
Theorem \ref{mt1}). It will be readily seen that this procedure is an extension
of the standard method of writing down an eigenbasis of $M(n)$ (see Chapter 2
in {\bf\cite{st}}).

Consider the decomposition
\beq \label{od} 
\C[B_q(n+1)] &=& (\C[B_q(n)]\oplus W(0)) \oplus 
\left(\oplus_{\chi \in \cN_q(n)} W(\chi)\right).
\eeq
We claim that
\beq 
\C[B_q(n)]\oplus W(0)\mbox{ and }W(\chi),\chi\in \cN_q(n) 
 \mbox{ are  $D_{n+1}$-closed}.
\eeq
This can be seen as follows. Consider the {\em standard inner product} on
$\C[B_q(n+1)]$ (i.e., declare $B_q(n+1)$ to be an orthonormal basis), which
is $GL(n+1, \F_q)$-invariant (and hence $H(n+1, \F_q)$-invariant). If follows that
$W(0)$, $W(\chi),\;\chi\in\cN_q(n)$ are orthogonal and hence the decomposition
(\ref{od}) is orthogonal.
Since $D_{n+1}$ is the adjoint of $U_{n+1}$, the claim now follows from
(\ref{d0}) and (\ref{d3}).

Thus $\C[B_q(n)]\oplus W(0)$ and $W(\chi),\chi\in\cN_q(n)$ are closed
under $\cA_q(n+1)$, which facilitates 
an inductive approach to the eigenvectors.

Let $X\in B_q(n,k)$ and consider $D_{n+1,k+1}(\theta_n(X))$. Then $\theta_n(X)$
is a sum of $q^{n-k}$ subspaces in $A_q(n+1,k+1)$ (from Lemma \ref{el}(iv)) 
and we can write
\beq \label{dtheta}
D_{n+1,k+1}(\theta_n(X)) = q^{n-k}X + v,
\eeq
where $v\in \C[A_q(n+1,k)]$. A little reflection shows that $v\in W(0)$.
Setting 
$$v=D^{'}_{n+1,k+1}(\theta_n(X))$$ 
gives a linear map
\beqn D^{'}_{n+1,k+1}: W(0)\rar W(0),\;\;0\leq k \leq n.
\eeqn

Define $\cA^{'}_q(n):W(0)\rar W(0)$ by
\beqn
\cA^{'}_q(n) = U_{n+1} + \sum_{k=0}^n q^k D^{'}_{n+1,k+1}.
\eeqn

We have the following relations (the first of which follows from
(\ref{dtheta})):
\beq \label{ind1}
\cA_q(n+1)(\theta_n(v))&=& q^n v 
       + \cA^{'}_q(n)(\theta_n(v)),\;\;v\in \C[B_q(n)],\\ \label{ind2}
\cA_q(n+1)(v)&=& \cA_q(n)(v) + \theta_n(v),\;\;v\in \C[B_q(n)].
\eeq

We now write down the matrix of $\cA^{'}_q(n)$ with respect to the basis
$\{\theta_n(X)\;|\;X\in B_q(n)\}$ of $W(0)$.

It follows from Lemma \ref{el}(iii) and Lemma \ref{orl}(ii) that
\beqn U_{n+1}(\theta_n(X)) &=& \sum_Y q\,\theta_n(Y),\;\;X\in B_q(n),
\eeqn
where the sum is over all $Y\in B_q(n)$ covering $X$. 
Similarly, it follows that
\beqn q^k D^{'}_{n+1,k+1}(\theta_n(Y))&=& 
\sum_X q \left\{ q^{k-1}\theta_n(X)\right\},\;\;Y\in B_q(n,k),
\eeqn
where the sum is over all $X\in B_q(n)$ covered by $Y$.

Thus we see that 
\beq \label{ind}
&\mbox{Matrix of $\cA^{'}_q(n)$ with respect to the basis
$\{\theta_n(X)\;|\;X\in B_q(n)\}$ is $qM_q(n)$.}
\eeq

Let $(V_1,f_1)$ be a pair consisting of a finite dimensional vector space
$V_1$ (over $\C$) and a linear operator $f_1$ on $V$. Let $(V_2,f_2)$ be
another such pair. By an isomorphism of pairs $(V_1,f_1)$ and $(V_2,f_2)$ 
we mean a linear isomorphism $\tau : V_1
\rar V_2$ such that $\tau(f_1(v)) = f_2(\tau(v)),\;v\in V_1$.

\bt \label{grind} 
Let $\chi\in\cN_q(n)$ and $X=X(\chi)$. Define
\beqn
\lambda(\chi): \C[B_q(X)]\rar W(\chi)
\eeqn
by $Y\mapsto q^{-\dim(Y)}p(\chi)(\wh{Y}),\;\;Y\in B_q(X)$.

Then 

(i) $\lambda(\chi)$ is an isomorphism of pairs $(\C[B_q(X)], qU_X)$
and $(W(\chi),U_{n+1})$.

(ii) $\lambda(\chi)$ is an isomorphism of pairs $(\C[B_q(X)], D_X)$
and $(W(\chi),D_{n+1})$.
\et
\pf
By Theorem \ref{f}(iii) it follows that $\lambda(\chi)(Y)\not= 0$ for all
$Y\in B_q(X)$. 
By (\ref{grd}) the dimensions of $\C[B_q(X)]$ and
$W(\chi)$ are the same. 
Thus, it follows from (\ref{ds}) that
$\lambda(\chi)$ is a vector space isomorphism.

(i) Let $Y\in B_q(X)$ with $\dim(Y)=k$.
 
We have (below the sum is over all $Z$ covering $Y$ in $B_q(X)$)
\beqn
\lambda(\chi)(qU_X(Y)) &=& 
q\lambda(\chi)\;\left(\sum_Z Z\right)\\
             &=& 
q^{-k}\sum_Z p(\chi)(\wh{Z}).
\eeqn

Before calculating $U_{n+1}\lambda(\chi)(Y)$ we make the
following observation. By Lemma \ref{el}(ii) every element covering $\wh{Y}$
is of the form $\wh{Z}$, for some $Z$ covering $Y$ in $B_q(n)$. Suppose
$Z\in B_q(n) - B_q(X)$. 
Since $\dim(W(\chi))= G_q(n-1)$, 
it follows by parts (ii) and (iii) of Theorem \ref{f} and (\ref{ds}) that
$p(\chi)(\wh{Z})=0$.

We now calculate $U_{n+1}\lambda(\chi)(Y)$. In the second step below we have
used the fact that $U_{n+1}$ is $H(n+1,\F_q)$-linear and in the third step,
using the observation in the paragraph above, we may restrict the sum to all
$Z$ covering $Y$ in $B_q(X)$.

We have
\beqn  
U_{n+1}(\lambda(\chi)(Y))&=& U_{n+1}\left( q^{-k}p(\chi)(\wh{Y}) \right)\\
    &=& q^{-k}p(\chi)(U_{n+1}(\wh{Y}))\\ 
    &=& q^{-k}\sum_Z p(\chi)(\wh{Z}).    
\eeqn

(ii) Let $Y\in B_q(X)$ with $\dim(Y)=k$.

We have (below the sum is over all $Z$ covered by  $Y$ in $B_q(X)$)
\beqn
\lambda(\chi)(D_X(Y)) &=& 
\lambda(\chi)\;\left(\sum_Z Z\right)\\
             &=& 
q^{-k+1}\sum_Z p(\chi)(\wh{Z}).
\eeqn

Before calculating $D_{n+1}\lambda(\chi)(Y)$ we make two observations:

(a) Let $Y$ cover $Z$, $Z\in B_q(X)$. Then, by Lemma \ref{orl}(ii),
there are $q$ subspaces in $[\wh{Z}]$ which are covered by $\wh{Y}$.
Let $Z_1\in [\wh{Z}]$ with $\wh{Y}$ covering $Z_1$. 
Then, there exists $g\in H(n+1,\F_q)$ with 
$g\wh{Z}=Z_1$.
It follows by Lemma \ref{el}(iii) that $g\wh{Y}=\wh{Y}$. Thus, from 
Lemma \ref{cl} we have $\chi(g)=1$.

(b) Let $\wh{Y}$ cover $Z$, where $Z\in B_q(n)$. Then $Z=Y$
and $p(\chi)(Z)=0$, since $\chi$
is nontrivial and every element of $H(n+1,\F_q)$ fixes $Z$.

Now we compute (using Lemma \ref{orl}(ii), Theorem \ref{f}(i), and (a), (b)
above)
\beqn  
D_{n+1}(\lambda(\chi)(Y))&=& 
q^{-k}\left\{D_{n+1}\left(p(\chi)(\wh{Y}) \right)\right\}\\
    &=& q^{-k}p(\chi)(D_{n+1}(\wh{Y}))\\ 
    &=& q^{-k+1}\sum_Z p(\chi)(\wh{Z}).    
\eeqn
where the sum is over all $Z\in B_q(X)$ covered by $Y$. \eprf

Before proceeding further we introduce some notation.
Let
$X\in B_q(n,n-1)$. The pairs $(\C[B_q(X)],U_X)$ and $(\C[B_q(n-1)],U_{n-1})$ 
are clearly isomorphic with many possible isomorphisms. We now define a
canonical isomorphism, based on the concept of a matrix in Schubert normal
form.

A $n\times k$ matrix $M$ over $\F_q$ is in {\em Schubert normal form}
 (or, {\em column reduced echelon form}) provided

\noi (i) Every column is nonzero.

\noi (ii) The last  nonzero entry in every column is a $1$. Let the last
nonzero entry in column $j$ occur in row $r_j$.

\noi (iii) We have $r_1 < r_2 < \cdots < r_k$ and the submatrix of $M$
formed by the rows $r_1,r_2,\ldots ,r_k$ is the $k\times k$ identity matrix.
We call $\{r_1,\ldots ,r_k\}$ the {\em pivotal indices} of $M$.

It is well known
that every $k$ dimensional subspace
of $\F^n_q$ is the column space of a unique $n\times k$ matrix in Schubert
normal form.
Given $X\in B_q(n,k)$, define $P(X)\seq \{1,2,\ldots ,n\}$ to be the pivotal
indices of the $n\times k$ matrix in Schubert normal form
with column space $X$. It is easy to see that $P(X)$
can also be defined
as follows
$$P(X)=\{j\in\{1,2,\ldots ,n\}\;:\; X\cap \F^j_q \in A_q(j)\}.$$
 
Let $X\in B_q(n,n-1)$ and let $M(X)$ be the $n\times (n-1)$ matrix in
Schubert normal form with column space $X$. The map $\tau(X) : \F^{n-1}_q
\rar X$ given by $e_j \mapsto \mbox{ column $j$ of }M(X)$ is clearly a linear
isomorphism and this isomorphism gives rise to an isomorphism
$$\mu(X) : \C[B_q(n-1)]\rar \C[B_q(X)]$$
of pairs $(\C[B_q(n-1)],U_{n-1})$ and $(\C[B_q(X)],U_X)$ 
(and also of pairs $(\C[B_q(n-1)],D_{n-1})$ and $(\C[B_q(X)],D_X)$) 
given by
$\mu(X)(Y)=\tau(X)(Y),\;Y\in B_q(n-1)$. It now follows from Theorem \ref{grind}
that

\bt \label{grind1}
Let $\chi\in \cN_q(n)$ and $X=X(\chi)$. Then the composition
$\lambda(\chi)\mu(X)$ is an isomorphism of pairs $(\C[B_q(n-1)],q\cA_q(n-1))$
and $(W(\chi),\cA_q(n+1))$. \eprf
\et

We shall now prove Theorems \ref{evectors} and \ref{mt1} by inductively writing
down eigenvectors of $\cA_q(n)$. We shall need an indexing set 
for the eigenvectors.
Given the fact that the multiplicities are the $q$-binomial coefficients 
it might appear that the set of subspaces $B_q(n)$ may be used as an indexing set.
We do not know of any natural way to index the eigenvectors of $\cA_q(n)$
by $B_q(n)$ (unlike the $q=1$ case, where the eigenvectors of $M(n)$ may be
naturally indexed by $B(n)$).  
A more useful 
indexing set, defined below, for the eigenvectors of $\cA_q(n)$
is suggested by the decomposition (\ref{od}).

For $n\geq 0$, inductively define a set $\cE_q(n)$ consisting of sequences
as follows (here $()$ denotes the empty sequence):
\beqn
\cE_q(0) &=& \{()\},\\
\cE_q(1) &=& \{(0),(1)\},\\
\cE_q(n) &=& 
\{(\alpha_1,\ldots ,\alpha_t)\;|\;
             (\alpha_1,\ldots \alpha_{t-1})\in
\cE_q(n-1),\;\alpha_t \in \{0,1\}\}\\
&&\cup \{(\alpha_1,\ldots ,\alpha_t)\;|\;
             (\alpha_1,\ldots \alpha_{t-1})\in\cE_q(n-2),\;
\alpha_t \in \cN_q(n-1)\},\;n\geq 2.
\eeqn
Given $\alpha\in \cE_q(n)$, let $N(\alpha)$ denote the number of nonzero
entries in the sequence $\alpha$ ( note that a nonzero entry is either 1 or
an element of $\cN_q(m)$ for some $m$). Set
\beqn
\cE_q(n,k) &=& \{\alpha\in \cE_q(n)\;|\;N(\alpha)=k\},\\
e_q(n,k) &=& |\cE_q(n,k)|.
\eeqn
It is easy to see that
\beq 
e_q(n+1,k) &=& e_q(n,k) + e_q(n,k-1) + (q^n - 1)e_q(n-1,k-1),\;n,k\geq
1,\eeq
with $e_q(0,k)=\delta(0,k)$ and $e_q(n,0)=1$,
the same recurrence (with the same initial conditions) as (\ref{nsw}), Thus
$e_q(n,k)=\qb{n}{k}$ and $|\cE_q(n)|=G_q(n)$.

We now prove Theorem \ref{mt1} in the following form, which gives extra
information in the form of eigenvectors.
\bt \label{mt}
For each $\alpha\in \cE_q(n)$ we define a vector $v_\alpha \in \C[B_q(n)]$
such that

(i) $\cA_q(n)(v_\alpha) = ((n-k)_q - (k)_q)v_\alpha$, where $k=N(\alpha)$.

(ii) $\{v_\alpha\;|\; \alpha\in \cE_q(n)\}$ is a basis of $\C[B_q(n)]$.
\et
\pf The proof is by induction on $n$, the cases $n=0,1$ being clear by
taking ($\{0\}$ denotes the zero subspace)
$$
v_{()}\;=\;\{0\},\;\;v_{(0)}\;=\;\{0\} + \F_q,\;\;v_{(1)}\;=\;\{0\} - \F_q.
$$

Let $n\geq 1$ and consider $\alpha=(\alpha_1,\ldots ,\alpha_t)\in \cE_q(n+1)$.
Set $\beta=(\alpha_1,\ldots ,\alpha_{t-1})$ and $k=N(\beta)$. We have three
cases:

(a) $\alpha_t=0$: We have $v_\beta\in \C[B_q(n)]$. Define
$$v_\alpha = q^kv_\beta + \theta_n(v_\beta) \in \C[B_q(n)]\oplus W(0).
$$

(b) $\alpha_t=1$: We have $v_\beta\in \C[B_q(n)]$. Define
$$v_\alpha = q^{n-k}v_\beta - \theta_n(v_\beta) \in \C[B_q(n)]\oplus W(0).
$$

(c) $\alpha_t = \chi,\;\chi\in\cN_q(n)$: We have $v_\beta\in \C[B_q(n-1)]$.
Set $X=X(\chi)$ and define
$$v_\alpha = \lambda(\chi)\mu(X)(v_\beta)\in W(\chi).
$$
Let us now check assertions (i) and (ii) in the statement of the theorem,
beginning with (i). We have three cases.

(a) $\alpha_t=0$: By the induction hypothesis and (\ref{ind}) we have
$$
\cA_q(n)(v_\beta)= ((n-k)_q - (k)_q)(v_\beta),
\;\;\;\cA^{'}_q(n)(\theta_n(v_\beta))= q((n-k)_q - (k)_q)(\theta_n(v_\beta)).
$$ 
We have, by (\ref{ind1}) and (\ref{ind2}),
\beqn
\cA_q(n+1)(v_\alpha)&=& \cA_q(n+1)(q^k v_\beta + \theta_n(v_\beta))\\
                    &=& q^k\cA_q(n+1)(v_\beta) 
                    + \cA_q(n+1)(\theta_n(v_\beta))\\
                    &=& q^k(\cA_q(n)(v_\beta) +\theta_n(v_\beta))
                        +q^n v_\beta + \cA^{'}_q(n)(\theta_n(v_\beta))\\
                    &=& (q^{n-k} + (n-k)_q - (k)_q)q^k v_\beta
                        + (q^k + q((n-k)_q - (k)_q))\theta_n(v_\beta)\\
                    &=& ((n+1-k)_q - (k)_q)(q^k v_\beta + \theta_n(v_\beta)).
\eeqn

(b) $\alpha_t=1$: By the induction hypothesis and (\ref{ind}) we have
$$
\cA_q(n)(v_\beta)= ((n-k)_q - (k)_q)(v_\beta),
\;\;\;\cA^{'}_q(n)(\theta_n(v_\beta))= q((n-k)_q - (k)_q)(\theta_n(v_\beta)).
$$ 
We have, by (\ref{ind1}) and (\ref{ind2}),
\beqn
\cA_q(n+1)(v_\alpha)&=& \cA_q(n+1)(q^{n-k}v_\beta - \theta_n(v_\beta))\\
                    &=& q^{n-k}\cA_q(n+1)(v_\beta) 
                    - \cA_q(n+1)(\theta_n(v_\beta))\\
                    &=& q^{n-k}(\cA_q(n)(v_\beta) +\theta_n(v_\beta))
                        - q^n v_\beta - \cA^{'}_q(n)(\theta_n(v_\beta))\\
                    &=& (-q^k + (n-k)_q - (k)_q)q^{n-k} v_\beta
                        - (-q^{n-k} + q((n-k)_q - (k)_q))\theta_n(v_\beta)\\
                    &=& ((n+1-(k+1))_q - (k+1)_q)(q^{n-k} v_\beta - \theta_n(v_\beta)).
\eeqn
(c) $\alpha_t = \chi,\;\chi\in\cN_q(n)$: Set $X=X(\chi)$. It follows from
Theorem \ref{grind1} that
\beqn \cA_q(n+1)(v_\alpha)&=&q((n-1-k)_q - (k)_q)v_\alpha \\
                           &=&((n+1-(k+1))_q - (k+1)_q)v_\alpha.
\eeqn

Assertion (ii) follows from the induction hypothesis using 
the decomposition (\ref{od}), the isomorphism
(\ref{d1}) and observing that the determinant of the $2\times 2$ matrix
$$ \left[ \ba{cr} q^k & 1 \\ q^{n-k} & -1 \ea \right]
$$
is nonzero. \eprf

We denote the basis given in part (ii) of Theorem \ref{mt} by $\cB_q(n)$. Note
that (upto scalars) this basis is canonical in the sense that we have not made
any choices anywhere.

\bt \label{mt3}
The basis $\cB_q(n)$ of $\C[B_q(n)]$ 
is orthogonal with respect to the inner product
(\ref{inp}). 
\et
\pf The proof is by induction on $n$, the cases $n=0,1$ being clear.

Let $n\geq 1$. We consider two cases:

(i) Let $\beta=(\beta_1,\ldots ,\beta_{t-1})\in\cE_q(n)$.
Set $k=N(\beta)$ and 
\beqn
&\alpha=(\beta_1,\ldots ,\beta_{t-1},0),\;\;\;\; 
\alpha'=(\beta_1,\ldots ,\beta_{t-1},1).
\eeqn

Given a vectors $u,v\in\C[B_q(n)]$, we shall write $\inp{u}{v}_n$ for
the inner product (\ref{inp}) calculated in $\C[B_q(n)]$ and $\inp{u}{v}_{n+1}$
for the inner product calculated in $\C[B_q(n+1)]$. We have, for
$X\in B_q(n,k)$,
\beqn 
&\inp{X}{X}_n \;=\; \frac{q^{\bin{k}{2}}}{P_q(n)},\;\; 
\inp{X}{X}_{n+1} \;=\; \frac{q^{\bin{k}{2}}}{P_q(n+1)}\;\;
                 \;=\;\frac{1}{1+q^n}\inp{X}{X}_n,&\\
&\inp{\theta_n(X)}{\theta_n(X)}_{n+1}\;=\;
\frac{q^{\bin{k+1}{2}}}{P_q(n+1)}\,q^{n-k}
\;=\;\frac{q^n}{1+q^n}\inp{X}{X}_n.&
\eeqn
It follows that
\beq \label{inpind}
&\inp{v}{v}_{n+1}\;=\;\frac{1}{1+q^n}\inp{v}{v}_n,\;\;
\inp{\theta_n(v)}{\theta_n(v)}_{n+1}\;=\;
\frac{q^n}{1+q^n}\inp{v}{v}_n,\;\;v\in \C[B_q(n)].&
\eeq
Note that the scalar factors on the right hand side are uniform across all vectors
and do not depend on $k$. Thus, since $\C[B_q(n)]$ and $W(0)$ are orthogonal
in $\C[B_q(n+1)]$, 
it follows by the induction hypothesis that
$\{v_\beta , \theta_n(v_\beta)\,|\, \beta\in \cE_q(n)\}$ is an orthogonal basis
of $\C[B_q(n)]\oplus W(0)$.

We have
\beq \label{indfs1}
&v_\alpha\;=\;q^kv_\beta + \theta_n(v_\beta),\;\;
v_{\alpha'}\;=\;q^{n-k}v_\beta - \theta_n(v_\beta).&
\eeq
Since $v_\beta$ is orthogonal to $\theta_n(v_\beta)$ we have,
using (\ref{inpind}),
\beqn
\inp{v_\alpha}{v_{\alpha'}}_{n+1} &= & 
q^n\inp{v_\beta}{v_\beta}_{n+1} - \inp{\theta_n(v_\beta)}
{\theta_n(v_\beta)}_{n+1}\\
&=&\frac{q^n}{1+q^n}\inp{v_\beta}{v_\beta}_n - 
\frac{q^n}{1+q^n}\inp{v_\beta}{v_\beta}_n\\
&=&0.
\eeqn
From the isomorphism $\theta_n$ we now see that
$$ \{ v_\alpha , v_{\alpha'}\;|\; \beta\in \cE_n(q)\}$$
is an orthogonal basis of $\C[B_q(n)]\oplus W(0)$.

(ii) Let $\beta = (\beta_1,\ldots ,\beta_{t-1})\in\cE_q(n-1)$ and 
let $\chi\in\cN_q(n)$.
Set $\alpha=(\beta_1,\ldots ,\beta_{t-1},\chi)\in\cE_q(n+1)$ 
and $X=X(\chi)$, where $X\in B_q(n,n-1)$. We have 
$v_\alpha = \lambda(\chi)\mu(X)(v_\beta)$.

Let $Y\in B_q(X)$ with $\dim(Y)=k$. We have
$$\inp{Y}{Y}_{n-1}\;=\;\frac{q^{\bin{k}{2}}}{P_q(n-1)}.
$$ 

Now observe the following: $p(\chi)(\wh{Y})$ is a linear combination of the elements of the orbit
$[\wh{Y}]$, whose cardinality is $q^{n-k}$. 
The number of elements $g\in H(n+1,\F_q)$ with $g\cdot \wh{Y}=\wh{Y}$
is $q^k$ and by Lemma \ref{cl} each such $g$ satisfies $\chi(g)=1$.
So, for $Z\in [\wh{Y}]$, if $g_1\cdot \wh{Y} = g_2 \cdot \wh{Y} = Z$ then
$\chi(g_1)=\chi(g_2)$.
Thus we have
\beqn
&\inp{\lambda(\chi)(Y)}{\lambda(\chi)(Y)}_{n+1}\;=\;
q^{-2k}\frac{q^{\bin{k+1}{2}}}{P_q(n+1)}\,q^{2k}q^{n-k}\;=\;
\frac{q^n}{(1+q^{n-1})(1+q^n)}\inp{Y}{Y}_{n-1}.&
\eeqn
It follows that
\beq \label{inpind1}
&\inp{\lambda(\chi)\mu(X)(v)}{\lambda(\chi)\mu(X)(v)}_{n+1}\;=\;
\frac{q^n}{(1+q^{n-1})(1+q^n)}\inp{v}{v}_{n-1},\;\;v\in \C[B_q(n-1)].&
\eeq
From the isomorphism $\lambda(\chi)\mu(X)$ we now see that
$$ \{ v_\alpha\;|\; \beta\in \cE_q(n-1)\}$$
is an orthogonal basis of $W(\chi)$.

That completes the proof. \eprf

The following result collects information about the length of the 
vectors $v_\alpha, \alpha
\in \cE_q(n)$, under the inner product (\ref{inp}) 
and the absolute values of their
standard coordinates $v_\alpha(Y), Y\in B_q(n)$ (i.e., 
$v_\alpha = \sum_Y v_\alpha(Y) Y$). 
Given $X\in A_q(n+1)$ we 
denote $X\cap \F^n_q$ by $X^r$. For $\alpha\in \cE_q(n)$, we denote by $\ol{\alpha}\in\cE_q(n)$ the sequence
obtained by interchanging the 0's and 1's in $\alpha$.

\bt \label{coeffl}

(a) Let $\beta= (\beta_1,\ldots ,\beta_{t-1})\in \cE_q(n)$ with $N(\beta)=k$.
We have 
 
(i) If $\alpha = (\beta_1,\ldots ,\beta_{t-1},0)$ then, for $Y\in B_q(n+1)$,
\beqn
v_\alpha (Y) &=& \left\{\ba{ll}
q^k v_\beta (Y)& \mbox{ if $Y\in B_q(n)$},\\
v_\beta (Y^r) & \mbox{ if $Y\in A_q(n+1)$}. \ea\right.
\\
\inp{v_\alpha}{v_\alpha}_{n+1}&=&\frac{q^n + q^{2k}}{1+q^n} \inp{v_\beta}{v_\beta}_n.
\eeqn

(ii) If $\alpha = (\beta_1,\ldots ,\beta_{t-1},1)$ then, for $Y\in B_q(n+1)$,
\beqn
v_\alpha (Y) &=& \left\{\ba{ll}
q^{n-k} v_\beta (X)& \mbox{ if $Y\in B_q(n)$},\\
-v_\beta (Y^r) & \mbox{ if $Y\in A_q(n+1)$}. \ea\right.
\\
\inp{v_\alpha}{v_\alpha}_{n+1}&=&\frac{q^n + q^{2(n-k)}}{1+q^n} \inp{v_\beta}{v_\beta}_n.
\eeqn

(b) Let $\beta = (\beta_1,\ldots ,\beta_{t-1})\in\cE_q(n-1)$ and 
let $\chi\in\cN_q(n)$. We have

If $\alpha=(\beta_1,\ldots ,\beta_{t-1},\chi)\in\cE_q(n+1)$ 
and $X=X(\chi)$, where $X\in B_q(n,n-1)$ then, for $Y\in B_q(n+1)$, 


\beqn
|v_\alpha (Y)| &=& \left\{\ba{ll} 0 & \mbox{ if $Y\in B_q(n)$},\\
|v_\beta (Y^r)|& \mbox{ if $Y\in A_q(n+1)$ and $Y^r \seq X$},\\
 0 & \mbox{ if $Y\in A_q(n+1)$ and $Y^r\not\seq X.$} \ea\right.
\\
\inp{v_\alpha}{v_\alpha}_{n+1}&=&\frac{q^n}{(1+q^{n-1})(1+q^n)} 
\inp{v_\beta}{v_\beta}_{n-1}.
\eeqn

(c) For $Y\in B_q(n)$ and $\alpha\in\cE_q(n)$ we have
$$|v_\alpha (Y)| = |v_{\,\ol{\alpha}} (Y)|,\;\;
\inp{v_\alpha}{v_\alpha}_{n+1}
= \inp{v_{\,\ol{\alpha}}}{v_{\,\ol{\alpha}}}_{n+1}.$$

\et
\pf Parts (a)(i) and (a)(ii) follow from (\ref{inpind}) and (\ref{indfs1}).

Now consider part (b). The formula for $\inp{v_\alpha}{v_\alpha}$ follows
from (\ref{inpind1}) and the formula for $|v_\alpha(Y)|$ follows from the
observation in the proof of case (ii) in Theorem \ref{mt3}. 

Part (c) follows easily by induction from parts (a) and (b) on
observing that  
$N(\ol{\alpha})= n- N(\alpha)$. 
\eprf 

We now single out a special set of $2^n$ eigenvectors of $M_q(n)$.
These can be seen as the $q$-analog of the classical eigenvectors
of $M(n)$, written down in {\bf\cite{cst,st}}. Define
$$\cI_q(n) = \{(\alpha_1,\ldots ,\alpha_n)\in \cE_q(n)\;:\;\alpha_i
\in \{0,1\} \mbox{ for all }i\}.
$$
Given $\alpha=(\alpha_1,\ldots ,\alpha_n)\in\cI_q(n)$ and $i\in\{1,2,\ldots ,n\}$,
define
$$d(\alpha,i)= |\{j<i\;:\;\alpha_j\not=\alpha_i\}|.$$
Note that, if $N(\alpha)=k$ then $\sum_{i=1}^n d(\alpha,i)=k(n-k)$.
Set $S(\alpha)=\{i\;:\;\alpha_i=1\}$.

Given $\alpha=(\alpha_1,\ldots ,\alpha_n)\in\cI_q(n)$ and $X\in B_q(n)$
set $d(\alpha,X)=\sum_i d(\alpha,i)$, where the sum is over all
$i\in \{1,2,\ldots ,n\}\setminus P(X)$.

\bt For $\alpha\in \cI_q(n)$ we have
\beqn v_\alpha &=& \sum_{X\in B_q(n)} (-1)^{|S(\alpha)\cap P(X)|}\;
q^{d(\alpha,X)}\;X.
\eeqn
In particular, $v_\alpha(\{0\})=q^{k(n-k)}$ and $v_\alpha(\F^n_q)=(-1)^k$,
where $k=N(\alpha)$.
\et

\noi {\bf{Remark}} Note that, when $q=1$, these are precisely the classical
eigenvectors of $M(n)$.

\pf By induction on $n$, the cases $n=0,1$ being clear. Let $n\geq 1$ and
consider $\alpha=(\alpha_1,\ldots ,\alpha_{n+1})\in \cI_q(n+1)$. Set
$\beta=(\alpha_1,\ldots ,\alpha_n)$ and consider $X\in B_q(n+1)$. 
We have the following cases:

(i) $n+1 \not\in P(X)$: We have $X=X\cap \F^n_q$. From the inductive hypothesis
and from cases (a), (b) in the proof of Theorem \ref{mt} we have
\beqn v_\alpha(X) &=& q^{d(\alpha,n+1)}\;v_\beta(X) \\
                  &=& (-1)^{|S(\beta)\cap P(X)|}\;q^{d(\alpha,n+1)}\;q^{d(\beta,X)}\\
                  &=& (-1)^{|S(\alpha)\cap P(X)|}\;q^{d(\alpha,X)}.
\eeqn

(ii) $n+1\in P(X)$ and $\alpha_{n+1}=1$: We have
\beqn v_\alpha(X) &=& - v_\beta(X^r) \\
                  &=& -(-1)^{|S(\beta)\cap P(X^r)|}\;q^{d(\beta,X^r)}\\
                  &=& (-1)^{|S(\alpha)\cap P(X)|}\;q^{d(\alpha,X)}.
\eeqn 

(iii) $n+1\in P(X)$ and $\alpha_{n+1}=0$: Similar to case (ii). \eprf

\bt The space of radial vectors $R_q(n)$ is contained in the subspace
spanned by $\{v_\alpha\;:\;\alpha\in \cI_q(n)\}$.
\et
\pf By induction on $n$, the cases $n=0,1$ being clear. Let $n\geq 1$.
By induction hypothesis and the isomorphism (\ref{d1}) we see that
\beqn
R_q(n)&\seq&\mbox{ Span }(\{v_\beta : \beta\in \cI_q(n)\}),\\
R_q(n+1)&\seq&\mbox{ Span }(\{v_\beta : \beta\in \cI_q(n)\})
\oplus \mbox{ Span }(\{\theta_n(v_\beta) : \beta\in \cI_q(n)\}),
\eeqn
and the right hand side of the second of these containments is equal to
$\mbox{ Span }(\{v_\alpha : \alpha\in\cI_q(n+1)\})$. \eprf

We are thus led to the following problem. For $0\leq k \leq n$, there is a unique
radial vector (up to scalars) that is an eigenvector of $\cA_q(n)$ with eigenvalue
$(n-k)_q - (k)_q$. Express this vector as a linear combination of
the vectors $\{v_\alpha : \alpha\in \cI_q(n), N(\alpha)=k\}$. The $n$-cube
case has a well known solution: the radial eigenvector is the sum of the vectors
in the (classical) eigenbasis with the same eigenvalue.

\section{{Weighted count of rooted spanning trees in $C_q(n)$}}

We now give the proof of the weighted count of rooted spanning trees
in $C_q(n)$. We use the definitions of Chapter 10 of {\bf\cite{st}}.

\noi
{\bf{Proof of Theorem \ref{qct}}} 
Form the directed loopless multigraph $D$ with $B_q(n)$ as the vertex set
and the following directed edges: for every edge $(X,Y)$ in $C_q(n)$ (where
we assume without loss of generality that $\dim(Y)=\dim(X)+1$) add $q^{\dim(X)}$
directed edges from $X$ to $Y$ in $D$ and one directed edge 
from $Y$ to $X$ in $D$. 
		
Now observe the following:

(i) The outdegree of a vertex $X$ in $D$ is 
$q^{\dim(X)}(n-\dim(X))_q + (\dim(X))_q = (n)_q$. Thus the matrix $L(D)$ (the
directed analog of the Laplacian) is given by
$$L(D)= (n)_qI - M_q(n).$$

(ii) There is an obvious root preserving onto map from the 
rooted oriented spanning subtrees of $D$ to the rooted 
spanning trees in $\cF_q(n)$,
where the inverse image of $F\in \cF_q(n)$ has cardinality $q^{w(F)}$.

(iii) By Theorem \ref{mt1}, the eigenvalues of $L(D)$ are 
$$(n)_q - ((n-k)_q - (k)_q) = (1+q^{n-k})(k)_q,\;\;k=0,1,\ldots ,n$$
with multiplicity $\qb{n}{k}$.

It follows from Theorem 10.4 in {\bf\cite{st}} (this is Tutte's directed
analog of the Matrix-Tree theorem) and item (ii) above
that the weighted count of rooted spanning trees in $\cF_q(n)$ is the product
of the nonzero eigenvalues of $L(D)$ and this agrees with the statement
of the Theorem by item (iii) above. \eprf 
\begin{center}
\section{{Acknowledgement}}
\end{center}

We are grateful to Professor Paul Terwilliger
for his encouragement and for detailed explanation on the origin of the
matrix $K_q(n)$. We thank Gaurav Bhatnagar for telling us about the paper
{\bf\cite{j}}. 

The first named author thanks the Indian Institute of Technology Bombay
for warm hospitality and support through the institute post-doctoral
fellowship.


\begin{thebibliography}{AAA}
%
%
\bibitem[A]{a} Askey, R.,
\newblock{\it Evaluation of Sylvester Type Determinants Using Orthogonal
Polynomials},
\newblock  in {\em Advances in Analysis}, World Scientific: 1-16  (2005).

\bibitem[CST]{cst} Ceccherini-Silberstein, T., Scarabotti, F.,
Tolli, F., {\em Harmonic analysis on finite groups, Representation theory, 
Gelfand pairs, and Markov chains},
Cambridge University Press (2008).
%
\bibitem[EK]{ek} Edelman, A., Kostlan, E.,
\newblock {\it How many zeros of a random polynomial are real?}, 
\newblock Bull. Amer. Math. Soc. (N.S.), {\bf 32} (1995), 1--37.


\bibitem[GR1]{gr1} Gasper, G., Rahman, M.,
\newblock Basic Hypergeometric Series,
\newblock{\it Cambridge University Press} (1990).

\bibitem[GS]{gs} Ghosh, S., Srinivasan, M. K.,
\newblock {\it A random walk on subspaces}, 
\newblock In preparation.


\bibitem[GR]{gr} Goldman, J., Rota, G. -C.,
\newblock{\it The number of subspaces of a vector space},
\newblock  in {\em Recent progress in Combinatorics} (Proc. Third Waterloo Conf.
on Combinatorics 1968), Academic Press : 75-83  (1969).



\bibitem[J]{j} Johnson, W. P.,
\newblock{ \it Some tridiagonal determinants},
\newblock The Ramanujan Journal, To appear.
Available at \url{https://doi.org/10.1007/s11139-021-00461-4}


\bibitem[K]{k} Kac, M.,
\newblock{ \it Random walk and the theory of Brownian motion},
\newblock Amer. Math. Monthly, 54: 369--391 (1947).

\bibitem[KC]{kc} Kac, V., Cheung, P.,
\newblock {\it Quantum Calculus},
\newblock  Springer-Verlag, 2002.


\bibitem[Ku]{ku} Kung, J. P. S.,
\newblock{\it The subset-subspace analogy},
\newblock  in {\em Gian-Carlo Rota on Combinatorics} 
(Contemporary mathematicians), Birkh\"{a}user Boston, : 277-283  (1995).



\bibitem[P]{p} Proctor, R. A.,
\newblock{ \it Representations of ${\mathfrak{sl}}(2,\C)$ on posets and the Sperner
property},
\newblock SIAM J. Alg. Discr. Methods, 3: 275--280 (1982).





\bibitem[S2]{s2} Srinivasan, M. K., 
\newblock{\it A positive combinatorial
formula for the complexity of the $q$-analog of the $n$-cube},
\newblock Electronic J. Comb., {\bf 19(2)} (2012), Paper 34 (14 Pages).

\bibitem[S3]{s3} Srinivasan, M. K., 
\newblock{\it The Goldman-Rota Identity and the Grassmann scheme},
\newblock Electronic J. Comb., {\bf 21(1)} (2014), Paper 37 (23 Pages).

\bibitem[St]{st} Stanley, R. P., 
\newblock Algebraic Combinatorics, Walks, Trees, Tableaux, and More (Second
Edition)
\newblock {\it Springer}, 2018.




\bibitem[TT]{tt} Taussky, O., Todd, J., 
\newblock {\it Another look at a matrix of Mark Kac}, 
\newblock Linear Algebra Appl., {\bf 150} (1991), 341--360.


\bibitem[T1]{t1} Terwilliger, P.,
\newblock{ \it The incidence algebra of a uniform poset},
\newblock in {\em Coding theory and design theory, Part I}, volume 20 of
{\em IMA Vol. Math. Appl.,} pages 193-212. Springer, New York, 1990.

\bibitem[T2]{t2} Terwilliger, P., 
\newblock{\it Two linear transformations each tridiagonal with respect to
an eigenbasis of the other; an algebraic approach to the Askey scheme of
orthogonal polynomials}, 
\newblock arXiv:math/0408390 (2004).

\bibitem[T3]{t3} Terwilliger, P., 
\newblock{\it Lowering-raising triples and $U_q(\mathfrak{sl}_2)$}, 
\newblock Linear Algebra Appl., {\bf 486} (2015), 1--172.

\bibitem[T4]{t4} Terwilliger, P., 
\newblock{\it Notes on the Leonard system classification}, 
\newblock Graphs Combin., {\bf 37} (2021), 1687--1748.



\end{thebibliography}
\end{document}